\documentclass[12pt]{article}
\title{Explicit Form of Coefficients in any MA(2) Process}
\author{Simon Ku and Eugene Seneta\footnote{Corresponding author. Email address: eseneta@maths.usyd.edu.au (Eugene Seneta).}\\{\small  School of Mathematics and Statistics FO7, University of Sydney, N.S.W. 2006, Australia}}

\begin{document}
\maketitle
\abstract
We shall show that for {\it any} $MA(2)$  process (apart from those with coefficients $\theta_1,\theta_2 $ lying on certain line-segments) there is {\it one and only one  invertible} $MA(2)$ process with the {\it same} autocovariances $\gamma_0,\gamma_1,\gamma_2$. 
It is this invertible version which computer-packages fit, regardless, even if data came from a non-invertible $MA(2)$ process.  This has consequences for prediction from a fitted process, inasmuch as such prediction would seem to be inappropriate.  We express the coefficients $\theta_1,\theta_2 $ of the invertible version in terms of  $\gamma_0,\gamma_1,\gamma_2$ explicitly using analytical reasoning, following a graphical approach of Sbrana (2012) which indicates this result within the invertibility region.  We also express $(\theta_1,\theta_2)$ in the non-invertibility region.
\section{$MA(2)$ Review}
Suppose we know that we are dealing with an $ MA(2)$ process:
 $$X(t)=e(t)-\theta_{1}e(t-1)-\theta_{2}e(t-2)$$
where $\sigma^2=Var(e(t))$, where $\{e(t)\}$ is white noise, $\theta_2 \neq 0,$ and whose non-zero autocovariances $\gamma_0,\gamma_1,\gamma_2$ are specified.

The ACF of $X(t)$ is 
\begin{eqnarray}
\frac{\gamma_1}{\gamma_0}&=&\rho_1=\frac{-\theta_1+\theta_1\theta_2}{1+\theta_1^2+\theta_2^2}  \nonumber \\
\frac{\gamma_2}{\gamma_0}&=&\rho_2=\frac{-\theta_2}{1+\theta_1^2+\theta_2^2}  \nonumber \\   
\frac{\sigma^2}{\gamma_0}&=&\frac{1}{1+\theta_1^2+\theta_2^2}.   \label{ma2_0}
\end{eqnarray}

Write
\begin{equation} \label{M1}
M(z)=z^2 - \theta_1z -\theta_2.
\end{equation} 
The $MA(2)$ process is {\it invertible} if and only if the roots $z_1, z_2$ (which are $\neq 0$ since $\theta_2 \neq 0$) of 
\begin{equation} \label{M2}
M(z) = 0 
\end{equation} 
satisfy $|z_i|<1, i =1,2.$

Equivalently, the invertibility conditions of $X(t)$, that is the region of $(\theta_1,\theta_2)$ in ${\bf R^2}$ which is commonly referred to as invertible triangle for $MA(2)$, are
\begin{eqnarray}
\theta_2-\theta_1&<&1  \label{ma2_00a}\\
\theta_2+\theta_1&<&1  \label{ma2_00b} \\   
\theta_2^2&<&1.   \label{ma2_00}
\end{eqnarray}

These can also be described as:
\begin{eqnarray} 
|1-\theta_2|&>& |\theta_1| \label{A}\\ 
|\theta_2|&<& 1. \label{B} 
\end{eqnarray}

In terms of $\gamma_1,\gamma_2$, and $\sigma^2$, $\theta_1$, and $\theta_2$ can be expressed as 

\begin{eqnarray}
\theta_1&=&-\frac{\gamma_1}{\sigma^2+\gamma_2}\nonumber \\
\theta_2&=&-\frac{\gamma_2}{\sigma^2}.          \label{ma2_1}
\end{eqnarray}

Give $\gamma_0,\gamma_1,\gamma_2$, the correct expression of $(\theta_1,\theta_2)$ depends on the correct choice of $\sigma^2$.

Substituting from (\ref{ma2_1}) in the definition of $\sigma^2$ in (\ref{ma2_0}), we have

\begin{eqnarray}
x^4+a_{1}x^3+a_{2}x^2+a_{1}kx+k^2=0 \label{ma2_2}
\end{eqnarray}
where $x=\sigma^2$ and $a_1=2\gamma_2-\gamma_0$, $a_2=2\gamma_2^2-2\gamma_0\gamma_2+\gamma_1^2$, $k=\gamma_2^2$, $a_3=a_{1}k$, $a_4=k^2$.

Sbrana (2011) (2012) asserts that there is only one solution of (\ref{ma2_2}), which he expresses explicitly in terms of $\gamma_0,\gamma_1,\gamma_2,$ which gives $(\theta_1,\theta_2)$ in (\ref{ma2_1}) corresponding to an invertible process (that is: satisfying (\ref{ma2_00a})-(\ref{ma2_00})). His reasoning is graphical (Sbrana, 2012), based on scanning Figure 3 of Stralkowski {\it et al.} (1974), which is Chart C, p.663 of Box {\it et al.} (2008). 

One motivation for the present paper is to verify this analytically. 

We shall first show that, for {\it any} $MA(2)$ process, apart from those with $(\theta_1,\theta_2)$ satisfying one of (a)-(c) below, there is {\it one and only one invertible} $MA(2)$ process with the {\it same} $\gamma_0,\gamma_1,\gamma_2$. 

It is this invertible version which computer-packages fit, regardless, even if data came from a non-invertible $MA(2)$ process.  This has consequences for prediction from a fitted process, inasmuch as such prediction would seem to be inappropriate.

We shall express $\theta_1,\theta_2,\sigma^2$ explicitly in terms of $\gamma_0,\gamma_1,\gamma_2$, for any $MA(2)$ process with an invertible version, irrespective of whether invertibility holds or not.

The (a)-(c) below correspond to there being a root of (\ref{M1}) on the unit circle: (a) a root $1$, (b) a root $-1$, and (c) a (complex) root $e^{i\lambda}, \lambda \neq 0, -\pi <\lambda < \pi$:

\begin{itemize}
\item[(a)] $1-\theta_1 - \theta_2=0.$
\item[(b)] $1+\theta_1 - \theta_2=0.$  
\item[(c)] $(\theta_1,\theta_2) = (2 cos \lambda, -1), \lambda \neq 0, -\pi <\lambda < \pi.$ 
\end{itemize}

Note that (a)-(c) together contain the {\it boundaries} of the invertibility triangle, the open set described by (\ref{ma2_00a})-(\ref{ma2_00}).

\section{ Anderson's Identity: Consequences}
Given an $MA(2)$ process with coefficients $\theta_1, \theta_2 \neq 0,\sigma^2,$ and autocovariances $ \gamma_0,\gamma_1,\gamma_2,$ a relationship between the two triples is given by Anderson's Identity (Anderson, 1971, Lemma 3.4.1):
\begin{equation} \label{A1}
\sum_{h=-2}^{2} \gamma_h z^h = \sigma^2M(z)M(z^{-1})
\end{equation} 
where $M(\cdot)$ is defined by (\ref{M1}). 

We now adapt to our specific situation the sketch-argument of Anderson (1971, Section 5.7). We were motivated by remarks in the paper of Ter\"asvirta (1977), within a more general setting.

The two roots of (\ref{M2}), $z_1,z_2$, are both non-zero, real, or are a complex conjugate pair. We may write:
$$M(z)M(z^{-1})=(z-z_1)(z-z_2)(z^{-1}-z_1)(z^{-1}-z_2)$$
so $z_1^{-1}, z_2^{-1}$ are the roots of $M(z^{-1})=0$.

Hence both $z_i$ and $z_i^{-1}$ are roots of 
\begin{equation} \label{A3}
\sum_{h=-2}^2 \gamma_hz^h = 0 
\end{equation} from (\ref{A1}). 
If $|z_i| \neq 1,i=1,2 $, then one of the roots $z_i,z_i^{-1}$, for each fixed $i$,  has absolute value {\it less than 1}.  Hence the four roots of (\ref{A3}) can be grouped into two sets $(w_1,w_2), (w_3,w_4)$, where $|w_i| <1,i=1,2$, and $|w_i|>1,i=3,4.$ 
Now define 
$$M^{*}(z)=(z-w_1)(z-w_2). $$
Then (\ref{A1}) holds, with $M(z)$ on the right-hand side replaced by $M^{*}(z)$ and $\sigma^2$ replaced by 
$$(\sigma^*)^2 = \sum_{h=-2}^2 \gamma_h /(M^*(1))^2. $$

Thus if $|z_i| \neq 1, i=1,2 $, that is: if there is no unit modulus root of $M(z)=0$, then there is an {\it invertible} $MA(2)$ process with these specified $\gamma_0,\gamma_1,\gamma_2$.

Pursuing the case $|z_i| \neq 1,i=1,2 $ further, we see that if $z_1$ is {\it real} and if $z_1 \neq z_2$, there are {\it  four distinct} $MA(2)$ processes with these same specified $\gamma_0, \gamma_1, \gamma_2$. These are defined by taking $M^*(z)= (z-v_1)(z-v_2)$, where 
$ (v_1,v_2) \in \{ (z_1,z_2), (z_1,z_2^{-1}), (z_2,z_1^{-1}), (z_1^{-1},z_2^{-1})\},$ with corresponding coefficients 
$$(\theta_1^*, \theta_2^*) = (v_1+v_2, -v_1v_2). $$

Note that we may not choose  $(v_1, v_2) = (z_i, z_i^{-1})$ to define $M^*(z)$, since then $M^*(z)M^*(z^{-1})$ would not involve $z_j, j \neq i $ at all, so (\ref{M1}) would not hold.

Next, if $|z_i| \neq 1, i=1,2 $ if $z_1$ is {\it real} and if $z_1 = z_2$, the above argument shows that there will be {\it just two distinct} $MA(2)$ processes with these same specified  $\gamma_0, \gamma_1,\gamma_2.$

If $z_1$ is {\it complex}, and $|z_i| \neq 1, i=1,2 $, then $z_1, z_2$ are complex conjugates, and there are {\it just two} $MA(2)$ processes each of form
with $ v_1=a^*e^{i\lambda^*}, \lambda^* \neq 0,  a^*=|v_1| \neq 1, v_2=a^*e^{-i\lambda^*}.$  The coefficients are, for each, of form:
$$(\theta_1^*, \theta_2^*) = (v_1+v_2, -v_1v_2) = (2a^* \cos \lambda^*, -(a^*)^2). $$

Finally, if $|z_i|=1$ for at least one of $i=1,2$, each of the possible choices of the pair $(v_1,v_2)$ to form $MA(2)$ processes with the prespecified $\gamma_0, \gamma_1, \gamma_2,$ will have $|v_i|=1$ for at least one of $i=1,2$.  Thus none of these processes will be invertible, and the coefficients $(\theta_1,\theta_2)$ of each are described by one of (a)-(c) above.  In particular, {\it there is no invertible version} if $|z_i|=1$ for at least one of $i=1,2.$

\section{Explicit Forms}
In this section we develop general theory, given any $\gamma_0,\gamma_1,\gamma_2$ for some $MA(2)$ process to express $\theta_1,\theta_2$ in terms of $\gamma_0,\gamma_1,\gamma_2$.  

Divided by $x^2$, the quartic equation (\ref{ma2_2}) is reduced to a quadratic equation in terms of $z$, where $z=x+\frac{k}{x}$,
\begin{eqnarray}
z^2+a_{1}z+(a_{2}-2k)=0 \label{ma2_3},
\end{eqnarray}
whence the roots of (\ref{ma2_3}) are
\begin{eqnarray}
z_{-}&=&\frac{1}{2}(-a_1-G)=\frac{1}{2}(\gamma_0-2\gamma_2-G)  \label{ma2_4a} \\
z_{+}&=&\frac{1}{2}(-a_1+G)=\frac{1}{2}(\gamma_0-2\gamma_2+G)        \label{ma2_4}
\end{eqnarray}
where

\begin{eqnarray}
G&=&\sqrt{a_1^2-4(a_2-2k)}\nonumber \\
&=&\sqrt{(2\gamma_2-\gamma_0)^2-4(\gamma_1^2-2\gamma_0\gamma_2) } \nonumber \\
&=&\sqrt{ 4\gamma_2^2+4\gamma_0\gamma_2+\gamma_0^2-4\gamma_1^2 } \nonumber \\
&=&\sqrt{(\gamma_0-2\gamma_1+2\gamma_2)(\gamma_0+2\gamma_1+2\gamma_2) }  \nonumber \\
&=&\sqrt{(\gamma_0-\frac{2(-\theta_1+\theta_1\theta_2)}{1+\theta_1^2+\theta_2^2}\gamma_0+2\frac{-\theta_2}{1+\theta_1^2+\theta_2^2}\gamma_0) }\times   \nonumber \\
& &\sqrt{(\gamma_0+\frac{2(-\theta_1+\theta_1\theta_2)}{1+\theta_1^2+\theta_2^2}\gamma_0+2\frac{-\theta_2}{1+\theta_1^2+\theta_2^2}\gamma_0)  } \nonumber \\
&=&\frac{\gamma_0}{1+\theta_1^2+\theta_2^2} \sqrt{(1-\theta_2+\theta_1)^2 (1-\theta_2-\theta_1)^2 }  \nonumber \\
&=&\frac{\gamma_0}{1+\theta_1^2+\theta_2^2}  |(1-\theta_2+\theta_1) (1-\theta_2-\theta_1) |  \nonumber \\
&=&\frac{\gamma_0}{1+\theta_1^2+\theta_2^2}  |(1-\theta_2)^2-\theta_1^2|.  \label{ma2_5}  
\end{eqnarray}

Note that $G^2$ is the discriminant of quadratic equation (\ref{ma2_3}), $G\ge 0$, and has been expressed in terms of $\gamma_i$, as well as in terms of $\theta_i$.

From $z=x+\frac{k}{x}$, we have
\begin{equation} \label{x2}
x^2-zx+k=0.
\end{equation}

Because $z$ can be taken $z_{-}$ or $z_{+}$, $x$ then can be four possible solutions, namely

\begin{eqnarray}
x_1&=&\frac{1}{2}(z_{-}-H_{-})     \label{ma2_8}\\
x_2&=&\frac{1}{2}(z_{-}+H_{-})     \label{ma2_9}\\
x_3&=&\frac{1}{2}(z_{+}-H_{+})     \label{ma2_10}\\
x_4&=&\frac{1}{2}(z_{+}+H_{+}).     \label{ma2_11}
\end{eqnarray}

Note that $x_1x_2=x_3x_4=k$, a property of a quadratic equation of (\ref{x2}). $H_{-}^2, H_{+}^2 $ is the discriminant of this quadratic equation in the respective cases $z = z_{-},z_{+}$.  Again, we shall show that $H_{-},H_{+}$ can be expressed in terms of the $\gamma_i$ as well as the $\theta_i$.

In general, given $x_i,i=1,2,3,4$, four sets of  $(\theta_1,\theta_2)$ can be defined as follows.

Taking $x_1$ as $\sigma^2$,
\begin{eqnarray}
\theta_1&=&-\frac{4\gamma_1}{\gamma_0+2\gamma_2-G-2H_{-}}     \nonumber \\
\theta_2&=&-\frac{4\gamma_2}{\gamma_0-2\gamma_2-G-2H_{-}}.     \label{ma2_22_1}
\end{eqnarray}  
 
Taking $x_2$ as $\sigma^2$,
\begin{eqnarray}
\theta_1&=&-\frac{4\gamma_1}{\gamma_0+2\gamma_2-G+2H_{-}}   \nonumber \\
\theta_2&=&-\frac{4\gamma_2}{\gamma_0-2\gamma_2-G+2H_{-}}.     \label{ma2_22_2}
 \end{eqnarray} 
 
Taking $x_3 $ as $\sigma^2$,
\begin{eqnarray}
\theta_1&=&-\frac{4\gamma_1}{\gamma_0+2\gamma_2+G-2H_{+}}    \nonumber \\
\theta_2&=&-\frac{4\gamma_2}{\gamma_0-2\gamma_2+G-2H_{+}}.     \label{ma2_22_3}
 \end{eqnarray} 
 
Taking $x_4$ as $\sigma^2$,
\begin{eqnarray}
\theta_1&=&-\frac{4\gamma_1}{\gamma_0+2\gamma_2+G+2H_{+}}     \nonumber \\
\theta_2&=&-\frac{4\gamma_2}{\gamma_0-2\gamma_2+G+2H_{+}}.     \label{ma2_22_4}
 \end{eqnarray} 

If $|1-\theta_2|>|\theta_1|$, that is for every $(\theta_1,\theta_2)$, {\it satisfying first two of invertibility conditions} (\ref{ma2_00a})-(\ref{ma2_00b}), then $G=\frac{\gamma_0}{1+\theta_1^2+\theta_2^2}  ((1-\theta_2)^2-\theta_1^2)$.

In terms of $\theta_1,\theta_2$, under $|1-\theta_2|>|\theta_1|$, from (\ref{ma2_4a})-(\ref{ma2_4}),

\begin{eqnarray}
z_{-}&=&\frac{1}{2}(\gamma_0-2\frac{-\theta_2}{1+\theta_1^2+\theta_2^2}\gamma_0-\frac{(1-\theta_2)^2-\theta_1^2}{1+\theta_1^2+\theta_2^2}\gamma_0) \nonumber \\
&=&\frac{\gamma_0}{1+\theta_1^2+\theta_2^2}(2\theta_2+\theta_1^2)        \label{ma2_6}
\end{eqnarray}
\begin{eqnarray}
z_{+}&=&\frac{1}{2}(\gamma_0-2\frac{-\theta_2}{1+\theta_1^2+\theta_2^2}\gamma_0+\frac{(1-\theta_2)^2-\theta_1^2}{1+\theta_1^2+\theta_2^2}\gamma_0) \nonumber \\
&=&\frac{\gamma_0}{1+\theta_1^2+\theta_2^2}(1+\theta_2^2)        \label{ma2_7}
\end{eqnarray}

and
\begin{eqnarray}
H_{-}&=&\sqrt{z_{-}^2-4k}=\sqrt{z_{-}^2-4\gamma_{2}^2}   \nonumber \\
&=&\sqrt{\frac{\gamma_0^2}{(1+\theta_1^2+\theta_2^2)^2}(2\theta_2+\theta_1^2)^2-4\frac{\theta_2^{2}\gamma_0^2}{(1+\theta_1^2+\theta_2^2)^2}} \nonumber \\
&=&\frac{\gamma_0}{1+\theta_1^2+\theta_2^2}\sqrt{\theta_1^2(\theta_1^2+4\theta_2)}        \label{ma2_12}
\end{eqnarray}

\begin{eqnarray}
H_{+}&=&\sqrt{z_{+}^2-4k}= \sqrt{z_{+}^2-4\gamma_{2}^2}   \nonumber \\
&=&\sqrt{\frac{\gamma_0^2}{(1+\theta_1^2+\theta_2^2)^2}(1+\theta_2^2)^2-4\frac{\theta_2^2\gamma_0^2}{(1+\theta_1^2+\theta_2^2)^2}} \nonumber \\
&=&\frac{\gamma_0}{1+\theta_1^2+\theta_2^2}\sqrt{(1-\theta_2^2)^2}    \nonumber \\
&=&\frac{\gamma_0}{1+\theta_1^2+\theta_2^2}|1-\theta_2^2|.  \label{ma2_13}
\end{eqnarray}

In view of (\ref{M1}), when the discriminant of the characteristic equation $I - \theta_1B -\theta_2B^2=0$ of $MA(2)$ satisfies:
\begin{eqnarray}
\theta_1^2+4\theta_2\geq 0 \label{ma2_14b}
\end{eqnarray}
then both $H_{-}$ of (\ref{ma2_12}) and $H_{+}$ of (\ref{ma2_13}) are real, and there are four real solutions, otherwise there are only two real solutions for $\sigma^2$.

If $|1-\theta_2|\le|\theta_1|$, then $G=\frac{\gamma_0}{1+\theta_1^2+\theta_2^2}  (\theta_1^2-(1-\theta_2)^2)$, and from (\ref{ma2_6})-(\ref{ma2_7}), $z_{-}=\frac{\gamma_0}{1+\theta_1^2+\theta_2^2}(1+\theta_2^2) $ and 
$z_{+}=\frac{\gamma_0}{1+\theta_1^2+\theta_2^2}(2\theta_2+\theta_1^2).$  Hence
$H_{-}=\frac{\gamma_0}{1+\theta_1^2+\theta_2^2}|1-\theta_2^2|$, $H_{+}=\frac{\gamma_0}{1+\theta_1^2+\theta_2^2}\sqrt{\theta_1^2(\theta_1^2+4\theta_2)} $.  Furthermore, when $|\theta_2|<1$ then $H_{-}=\frac{\gamma_0}{1+\theta_1^2+\theta_2^2}(1-\theta_2^2)$, otherwise $H_{-}=\frac{\gamma_0}{1+\theta_1^2+\theta_2^2}(\theta_2^2-1)$ when $|\theta_2|\ge 1$.

\section{Invertible $MA(2)$}

Suppose the process is known to invertible, so (\ref{ma2_00a})-(\ref{ma2_00}) all hold, and additionally that (\ref{ma2_14b}) holds.  Then (since $\theta_2 \neq 0$) $x_i,i=1,2,3,4$ are all positive:


\begin{eqnarray}
x_1&=&\frac{\gamma_0}{2(1+\theta_1^2+\theta_2^2)}(2\theta_2+\theta_1^2-\sqrt{\theta_1^2(\theta_1^2+4\theta_2)})  \nonumber \\     &=&\frac{\gamma_0}{1+\theta_1^2+\theta_2^2}{\bigg(}\sqrt{\frac{\theta_1^2}{4}}-\sqrt{\frac{\theta_1^2}{4}+\theta_2}\,\,{\bigg)}^2  \nonumber \\    
  &>&0    \label{ma2_15} \\
x_2&=&\frac{\gamma_0}{2(1+\theta_1^2+\theta_2^2)}(2\theta_2+\theta_1^2+\sqrt{\theta_1^2(\theta_1^2+4\theta_2)}) \nonumber \\   &=&\frac{\gamma_0}{1+\theta_1^2+\theta_2^2}\bigg{(} \sqrt{\frac{\theta_1^2}{4}}+\sqrt{\frac{\theta_1^2}{4}+\theta_2}\,\,{\bigg)}^2  \nonumber \\    
  &>&0    \label{ma2_16}\\
x_3&=&\frac{\gamma_0}{2(1+\theta_1^2+\theta_2^2)}(1+\theta_2^2-1+\theta_2^2) \nonumber \\   
&=&\frac{\gamma_0}{1+\theta_1^2+\theta_2^2}\theta_2^2  \nonumber \\    
  &>&0    \label{ma2_17}\\
x_4&=&\frac{\gamma_0}{2(1+\theta_1^2+\theta_2^2)}(1+\theta_2^2+1-\theta_2^2) \nonumber \\   
&=&\frac{\gamma_0}{1+\theta_1^2+\theta_2^2} \nonumber \\      
  &>&0.   \label{ma2_18}
\end{eqnarray}

{\it Under the invertibility conditions alone (i.e. irrespective of whether (\ref{ma2_14b}) holds) we see from the above that  $0< x_3 < x_4 $, and $x_4=\sigma^2$ is the only correct solution}.

If and only if additionally (\ref{ma2_14b}) holds, $x_1$ and $x_2$ are both real, and, in the event clearly $0 < x_1 \leq  x_2$.  The inequality is strict if the inequality in (\ref{ma2_14b}) is strict, as we shall assume for the rest of this section and (for convenience)  \S 5.



In fact then  $x_4=\max(x_1,x_2,x_3,x_4)$, $x_3=\min(x_1,x_2,x_3,x_4)$ and $x_3<x_1<x_2<x_4$. 
Since $\theta_2<1+\theta_1$,  $\theta_2<1-\theta_1$, therefore $\theta_2<1-|\theta_1|=1-2\frac{|\theta_1|}{2}$.  This leads to

\begin{eqnarray}
x_2&=&\frac{\gamma_0}{1+\theta_1^2+\theta_2^2}{\Bigg(}\sqrt{\frac{\theta_1^2}{4}}+\sqrt{\frac{\theta_1^2}{4}+\theta_2}\,{\Bigg)}^2  \nonumber \\
&<&\frac{\gamma_0}{1+\theta_1^2+\theta_2^2}{\Bigg(}\sqrt{\frac{\theta_1^2}{4}}+\sqrt{\frac{\theta_1^2}{4}+(1-2|\frac{\theta_1}{2}|)} \,{\Bigg)}^2  \nonumber \\
&=&\frac{\gamma_0}{1+\theta_1^2+\theta_2^2}{\Bigg(}\sqrt{\frac{\theta_1^2}{4}}+\sqrt{(1-|\frac{\theta_1}{2}|)^2}\, {\Bigg)}^2  \nonumber \\
&=&\frac{\gamma_0}{1+\theta_1^2+\theta_2^2} {\Bigg(}|\frac{\theta_1}{2}|+1- |\frac{\theta_1}{2}|\,{\Bigg)}^2  \nonumber \\
&=&x_4 \nonumber
\end{eqnarray}
and
\begin{eqnarray}
x_1&=&\frac{\gamma_0}{1+\theta_1^2+\theta_2^2}{\Bigg(}\sqrt{\frac{\theta_1^2}{4}}-\sqrt{\frac{\theta_1^2}{4}+\theta_2}\, {\Bigg)}^2 \nonumber \\
&>&\frac{\gamma_0}{1+\theta_1^2+\theta_2^2}{\Bigg(}\sqrt{\frac{\theta_1^2}{4}}-(\sqrt{\frac{\theta_1^2}{4}}+\sqrt{|\theta_2|})\,{\Bigg)}^2 \nonumber \\
&=&\frac{\gamma_0}{1+\theta_1^2+\theta_2^2}|\theta_2| \nonumber \\
&>&\frac{\gamma_0}{1+\theta_1^2+\theta_2^2}\theta_2^2 \nonumber \\
&=&x_3     \nonumber  
\end{eqnarray}
as $|\theta_2|<1$.  Here inequality of $\sqrt{a+b}<\sqrt{a}+\sqrt{|b|}$, if $a >0$, is used.

The ranking $x_3<x_1<x_2<x_4$ is, naturally, consistent with $x_{1}x_{2}=x_{3}x_{4}=k$.  Thus if $x_4$ is the largest of the four positive numbers, $x_3$ must be the smallest by this identity, and $x_1,x_2$ must be in between them in size.  This order holds for {\it any} $MA(2)$ processes (invertible or non-invertible).

In terms of $\gamma_0,\gamma_1,\gamma_2$ of an {\it  invertible}  $MA(2$) process, $(\theta_1, \theta_2)$ (without explicit involvement of $\sigma^2=x_4$ and irrespective of whether (\ref{ma2_14b}) holds or not) is therefore defined by (\ref{ma2_22_4}), 
where $G, H_{+}$ are also functions of $\gamma_i,i=0,1,2$.  That is to say $(\theta_1,\theta_2)$ constructed as above via $\sigma^2=x_4$ lie {\it always} in the invertible triangle.

\section{Identification of $\sigma^2$}   
Given $\gamma_0,\gamma_1,\gamma_2$,$x_i,i=1,2,3,4$ are defined by the general expressions (\ref{ma2_8})-(\ref{ma2_11}).  Label (\ref{ma2_00a})-(\ref{ma2_00}) as $(A),(B),(C)$.

If an $MA(2)$ process has an invertible version, then, according to our \S 2 there are eight cases for the pair $(\theta_1,\theta_2)$ which need to be checked to see for each case when $\sigma^2$ is calculated correctly.  In the case of invertibility, as we have seen the correct $\sigma^2$ is the maximum in magnitude (that is, of rank 4) of four possibilities $x_i,i=1,2,3,4$, as in (\ref{ma2_15})-(\ref{ma2_18}), if all four are real (and positive).  But in other cases, the size of correct $\sigma^2$ in magnitude can be of rank 1, 2, or 3.

By the superscript $c$ we shall mean the {\it strict} reverse inequality.  Thus while $A$ means $ \theta_2 - \theta_1 <1$,  $A^c$ will mean $ \theta_2 - \theta_1 >1$.  We then have the following eight cases: 

 \begin{eqnarray}
Case (1)&=&ABC  \nonumber \\ 
Case (2)&=&A^{c}BC  \nonumber \\ 
Case (3)&=&AB^{c}C  \nonumber \\ 
Case (4)&=&A^{c}B^{c}C  \nonumber \\ 
Case (5)&=&ABC^{c}  \nonumber \\ 
Case (6)&=&A^{c}BC^{c}  \nonumber \\ 
Case (7)&=&A^{c}B^{c}C^{c}  \nonumber \\ 
Case (8)&=&AB^{c}C^{c}  \label{ma2_23} 
\end{eqnarray}

Note that Case(4) is impossible.  If $D$ is defined to be (\ref{ma2_14b}) with strict inequality, then Case(1) and Case(5) can each be split, into  Case(1a) and Case(5a), that is $ABCD$ and $ABC^{c}D$ respectively; and Case(1b) and Case(5b), that is $ABCD^{c}$ and $ABC^{c}D^{c}$ respectively.  Only Case(1) and Case(5) may be so split, of the possible seven cases.  For the other five possible cases, $D$ must hold automatically.

The order of $0<x_3<x_1<x_2<x_4$ (or $0 <x_3< x_4$, if only $ x_3,x_4$ are real) holds even in the cases other than invertibility. 

Only two real solutions $x_3,x_4$ for $\sigma^2$ are obtained when $D^c$ holds.

Individual algebraic consideration of the above cases via $x_i,i=1,2,3,4$, gives the correct $\sigma^2$ as:
in Case(1a),  $x_4$, the largest of four; 
in Case(1b),  $x_4$, the largest of two ($x_3, x_4$);
in Case(2), $x_2$, the second largest of four; 
in Case(3), $x_2$, the second largest of four;  
in Case(5a),  $x_3$, the smallest of four; 
in Case(5b),  $x_3$, the smallest of two ($x_3, x_4$); 
in Case(6), $x_1$, the third largest of four; 
in Case(7),  $x_3$, the smallest of four; 
in Case(8),  $x_1$, the third largest of four.

Thus, given $\gamma_0,\gamma_1,\gamma_2$, and a case number for $(\theta_1,\theta_2)$, there are multiple candidates for $x_i$, but only one is ``correct" for $\sigma^2$.

We can actually simplify: only four situations are needed to classify all cases based on the constraints of $(\theta_1,\theta_2)$ in view of (\ref{A})-(\ref{B}).

If $|1-\theta_2|>|\theta_1|, |\theta_2|<1$, i.e. Case(1), that is Case(1a) and Case(1b), then use $x_4$ as the correct $\sigma^2$, the largest of $x_i$.

If $|1-\theta_2|>|\theta_1|, |\theta_2|\ge 1$, i.e. Case(5), that is Case(5a) and Case(5b), or Case(7), then use $x_3$ as the correct $\sigma^2$, the smallest of $x_i$ (the sign of $H_{+}$ in $z_{+}$ exchanged, so the correct $\sigma^2$ changes from $x_4$ (in Case(1)) to $x_3$).

If $|1-\theta_2|<|\theta_1|, |\theta_2|<1$, i.e. Case(2) or Case(3), then use $x_2$ as the correct $\sigma^2$, the second largest of $x_i$ (the sign of $H_{-}$ in $z_{-}$ exchanged, hence $z_{-}, z_{+}$ exchanged, so the correct $\sigma^2$ changes from $x_4$ (in Case(1)) to $x_2$).

If $|1-\theta_2|<|\theta_1|, |\theta_2| > 1$, i.e. Case(6) or Case(8), then use $x_1$ as the correct $\sigma^2$, the third largest of $x_i$ (the sign of $H_{+}$ in $z_{+}$, and  $H_{-}$ in $z_{-}$ both exchanged, so the correct $\sigma^2$ changes from $x_4$ (in Case(1)) to $x_1$).

We mention an especially interesting case of a non-invertible $MA(2)$ process which has an invertible version as indeed foreshadowed in our \S2: $\theta_2 = -1, |\theta_1| > 2.$  In this situation $H_ -=0,$ and $x_1=x_2$ is the ``correct"  $\sigma^2,$ and $x_3 < x_1 = x_2 < x_4$.

We see from the above that in any {\it non-invertible} case for which there is an invertible version, the largest $x_i$, namely $x_4$, is {\it never} the ``correct" $\sigma^2.$

For any such case, now choose an $x_i,i=1,2,3$ which is not the ``correct" $ x_i$  for the process, put it equal to $\sigma^2$ and construct $(\theta_1, \theta_2) $ via (\ref{ma2_1}).  If the resulting process were invertible, the ``correct" $\sigma^2$ would be $x_4$, by our \S4, a contradiction to our choice of $x_i$.  

We know that there is an invertible version from our \S 2, so, by elimination of $x_1, x_2, x_3$, it must correspond to $x_4$.
Now put $x_4 = \sigma^2$ and construct $ (\theta_1, \theta_2) $ via (\ref{ma2_1}), to give the invertible version of the given $MA(2)$ process with the same $\gamma_0, \gamma_1, \gamma_2.$

\section{Conclusions}

\begin{enumerate}
\item For the first time all solutions $x=\sigma^2$ of the quartic equation (\ref{ma2_2}) are considered, and each of them represents an $MA(2)$ invertible or non-invertible process, explicitly in terms of $(\theta_1,\theta_2)$ by (\ref{ma2_22_1})-(\ref{ma2_22_4}).
\item Given $\gamma_0,\gamma_1,\gamma_2$ corresponding to some $MA(2)$ process, there are either two or four $MA(2)$ processes with these autocovariances, precisely one of which is invertible.
\item Given $\gamma_0,\gamma_1,\gamma_2$ corresponding to some $MA(2)$ process, the unique invertible $MA(2)$ process with this autocovariance structure has $\sigma^2 =x_4$, where $x_4$ is explicitly given in terms of $\gamma_0,\gamma_1,\gamma_2$ by (\ref{ma2_11}), and the corresponding $(\theta_1,\theta_2)$ by (\ref{ma2_22_4}).
\item Given $\gamma_0,\gamma_1,\gamma_2$ corresponding to some $MA(2)$ process, providing, in addition, we know to which of the above seven possible cases in (\ref{ma2_23}) it corresponds, all of $\sigma^2,\theta_1,\theta_2$ can be explicitly specified and uniquely determined. 
\end{enumerate}

\section{Acknowledgement}
We thank Giacomo Sbrana for correspondence (24 October-13 November, 2012)  relating to the papers Sbrana (2011) (2012).  Our work on related topics began with Ku and Seneta (1998).

\end{document}